\documentclass{amsart}

\usepackage{amsmath}
\usepackage{amscd}
\usepackage{amssymb} 

\newcommand{\cal}{\mathcal}
\newcommand{\bk}{{\bf k}}

\newcommand{\bt}{{\bf t}}

\newcommand{\bQ}{{\Bbb Q}}

\newcommand{\bZ}{{\Bbb Z}}
\newcommand{\cA}{{\cal A}}

\newcommand{\cG}{{\cal G}}

\newcommand{\DGBV}{{\cal D}{\cal G}{\cal B}{\cal V}}

\newcommand{\fm}{{\frak m}}

\DeclareMathOperator{\ad}{ad}

\DeclareMathOperator{\Hom}{Hom}
\DeclareMathOperator{\Img}{Im}

\DeclareMathOperator{\Ker}{Ker}

\newtheorem{theorem}{Theorem}[section]
\newtheorem{proposition}{Proposition}[section]
\newtheorem{lemma}{Lemma}[section]

\theoremstyle{remark}

\theoremstyle{definition}
 
\newtheorem{definition}{Definition}[section]

\begin{document}
\title{On quasi-isomorphic DGBV algebras}
\author{Huai-Dong Cao and Jian Zhou}
\address{Department of Mathematics\\
Texas A\&M University\\
College Station, TX 77843}
\email{cao@math.tamu.edu, zhou@math.tamu.edu}
\begin{abstract}
One of the methods to obtain Frobenius manifold structures is via DGBV (differential 
Gerstenhaber-Batalin-Vilkovisky) algebra construction. An important problem is 
how to identify Frobenius manifold 
structures constructed from two different DGBV algebras. For DGBV algebras with suitable conditions,
we show the functorial property of a construction of deformations of
the multiplicative structures of their cohomology. In particular, we show that
quasi-isomorphic DGBV algebras yield identifiable Frobenius manifold 
structures.

\end{abstract}
\maketitle

String theorists are interested in two kinds of  
conformal field theories defined on a Calabi-Yau manifold $X$:
an $A$-type theory which depends only on the K\"{a}hler structure 
but not the complex structure $X$,
and an $B$-type theory which depends only on the complex structure 
but not the  K\"{a}hler structure.
Conceivably,
an $A$-type theory should then be related to the deformations of the
K\"{a}hler structure,
while an $B$-type theory should be related to 
the deformations of the complex structure.
The mysterious mirror symmetry \cite{Yau} can be formulated as 
the identification of an $A$-type theory
on $X$ with a $B$-type theory on its mirror manifold $\widehat{X}$.

Physicists also provide us some examples of such theories:
 the topological sigma model \cite{Wit} 
and the K\"{a}hler theory of gravity \cite{Ber-Sad} are $A$ type theories,
while the Kodaira-Spencer theory of gravity \cite{Ber-Cec-Oog-Vaf}
is a $B$ type theory.
Through the efforts of many mathematicians,
the topological sigma model now has rigorous mathematical formulation  
in terms of suitably defined Gromov-Witten invariants 
and has led to vast progress in symplectic geometry and algebraic geometry.
On the other hand, based on the work of Tian \cite{Tia} and Todorov \cite{Tod},
the Kodaira-Spencer theory of gravity was analyzed in
details by Bershadsky-Cecotti-Ooguri-Vafa \cite{Ber-Cec-Oog-Vaf},
and the theory of K\"{a}hler gravity by Bershadsky and Sadov \cite{Ber-Sad}.
Barannikov-Kontsevich \cite{Bar-Kon} reformulated the results in 
\cite{Ber-Cec-Oog-Vaf} in terms of Frobenius manifolds introduced by 
Dubrovin \cite{Dub1, Dub2}, and made the important observation  that 
there is an algebraic structure called DGBV algebra hidden in the theory,
and the method to obtain formal Frobenius manifold structure by 
the Kodaira-Spencer Lagrangian can be generalized 
to any DGBV algebra satisfying certain conditions.
See the detailed account in Manin \cite{Man2}.
In two earlier papers \cite{Cao-Zho1, Cao-Zho2},
we pointed out two DGBV algebra
structures in the theory of K\"{a}hler gravity, one on Dolbeault cohomology
and the other on de Rham cohomology,
and showed that they satisfy the conditions in \cite{Bar-Kon, Man2}
for constructing Frobenius manifold structures.
Furthermore,
we were able to identify the Frobenius manifold structures 
from these two different DGBV algebras.
Subsequently, we also generalized these results to hyperk\"{a}hler manifolds 
\cite{Cao-Zho3} and equivariant cohomology \cite{Cao-Zho4}.

We conjectured that for a Calabi-Yau manifold $X$
with a mirror manifold $\widehat{X}$,
one should be able to identify the Frobenius manifold structure constructed
in \cite{Bar-Kon} for $X$ with that constructed in \cite{Cao-Zho1, Cao-Zho2} 
for $\widehat{X}$ (maybe after some coordinate change). A general question
in this direction is how one can possibly identify Frobenius manifold structures constructed from two 
different DGBV algebras. Our idea is to first give a natural definition of 
quasi-isomorphisms of DGBV algebras, and then show that, when the constructions 
of formal Frobenius manifolds are applicable,
quasi-isomorphic DGBV algebras yield identifiable Frobenius manifold structures.
These are carried out in this paper (see {\em Definition 3.1} and 
{\em Theorem 3.2}).
We leave the problem of showing that 
the relevant DGBV algebras on $X$ and $\widehat{X}$ are quasi-isomorphic
to future investigation.

In the course of our study, we find it natural to consider homomorphisms of DGBV 
algebras.
The construction of Frobenius manifold structures is not functorial
with respect to general DGBV algebra homomorphisms,
but the construction of some one-parameter formal deformations is.
We establish the functorial properties by
showing the gauge invariance of the constructions.
Note that the gauge invariance has been studied in the special cases
of Kodaira-Spencer gravity \cite{Ber-Cec-Oog-Vaf} and K\"{a}hler gravity \cite{Ber-Sad}.

\section{Deformations of cohomology algebras of DGBV algebras}
\label{sec:deform}

Let $\bk$ be a graded commuative associative algebra with unit over $\bQ$,
$(\cA, \wedge)$ a graded 
commutative associative algebra with unit $1$ over ${\bk}$. 
For any linear operator $\Delta$ of odd degree,
define
\begin{eqnarray} \label{eqn:bracket}
[a \bullet b]_{\Delta} = (-1)^{|a|} (\Delta (a \wedge b) 
- (\Delta a) \wedge b - (-1)^{|a|} a \wedge \Delta b),
\end{eqnarray}
for homogeneous elements $a, b \in \cA$.
When there is no risk of confusion,
we will simply use $[\cdot \bullet \cdot]$ 
to denote $[\cdot\bullet\cdot]_{\Delta}$.
If $\Delta^2 = 0$ and
\begin{eqnarray} \label{eqn:derivation}
[a \bullet (b \wedge c)] 
	= [a \bullet b] \wedge c
	+ (-1)^{(|a|-1)|b|} b \wedge [a \bullet c],
\end{eqnarray}
for all homogeneous $a, b, c \in \cA$,
then $(\cA, \wedge, \Delta, [\cdot\bullet\cdot])$ is 
a {\em Gerstenhaber-Batalin-Vilkovisky (GBV) algebra}.
Under the above conditions,
it is straightforward to see that
\begin{eqnarray*}
&& [a \bullet b] = - (-1)^{(|a|-1)(|b|-1)}[b \bullet a], \\
&& [a \bullet [b \bullet c]] = [[a \bullet b] \bullet c]
+ (-1)^{(|a| -1)(|b| -1)} [b \bullet [a \bullet c]], \\
&& \Delta [a\bullet b] 
= [\Delta a \bullet b] + (-1)^{|a|+1} [a \bullet \Delta b].
\end{eqnarray*}
A  {\em DGBV (differential Gerstenhaber-Batalin-Vilkovisky) algebra}
is a GBV algebra with a $\bk$-linear derivation $\delta$ of odd degree
with respect to $\wedge$, such that
$$ \delta^2 = \delta\Delta + \Delta \delta = 0.$$
It is easy to see that
\begin{eqnarray*}
&& \delta [a\bullet b] 
= [\delta a \bullet b] + (-1)^{|a|+1} [a \bullet \delta b].
\end{eqnarray*}
It is routine to define homomorphisms of DGBV algebras.
Denote by $\DGBV$ the category of DGBV algebra.
One can define direct sum and tensor product in this category.

Let $\delta_{x(t)} = \delta + [x(t) \bullet \cdot]$,
where $t$ is an indeterminate and 
$x(t) \in \cA[[t]]$ is a formal power series in $t$ with coefficients
even elements of $\cA$.
Because of (\ref{eqn:derivation}), 
$\delta_x$ is a derivation of $\cA[[t]]$.
Now
$$
\delta_{x(t)}^2 
= [(\delta x(t) + \frac{1}{2} [x(t) \bullet x(t)]) \bullet \cdot].
$$
Hence if
\begin{eqnarray} \label{eqn:MC}
\delta x(t) + \frac{1}{2} [x(t) \bullet x(t)] = 0,
\end{eqnarray}
then $\delta_{x(t)}$ is a differential.
It is easy to see that if $\Delta x(t) = 0$, 
then $\delta_{x(t)} \Delta = - \Delta \delta_{x(t)}$.
As a consequence, 
$(\cA[[t]], \wedge, \delta_{x(t)}, \Delta, [\cdot \bullet \cdot])$ 
is a DGBV algebra.

We will be concerned with the deformations of the multiplicative structure
on the cohomology $H = H(\cA, \delta)$.
Since $\delta$ is a derivation, i.e.,
$$\delta (a \wedge b) = \delta a \wedge b + (-1)^{|a|} a \wedge \delta b,$$
for homogeneous $a, b \in \cA$,
$\wedge$ induces an associative product (also denoted by $\wedge$) on $H$.
Clearly, $(H, \wedge)$ is graded commutative.
Since we have  $\delta 1 = 0$,
the class of $1$ gives a unit of $(H, \wedge)$.
To summarize, $(H, \wedge)$ is a graded commutative associative algebra 
with unit over $\bk$.
Note that $(\cA[[t]], \wedge, \delta_{x(t)}, \Delta, [\cdot \bullet \cdot])$
is a formal deformation of 
$(\cA, \wedge, \delta, \Delta, [\cdot \bullet \cdot])$.
An important idea in \cite{Bar-Kon} and \cite{Man2}
is to obtain formal deformations of $(H, \wedge)$
by considering $H(\cA[[t]], \delta_{x(t)})$.
We will need the following

\begin{lemma} \label{lm:quasi}
Let $\cA$ be a vector space with two endomorphisms 
$\delta$ and $\Delta$ satisfying 
$\delta^2 = \Delta^2 = \delta \Delta + \Delta \delta = 0$.
Then $\Img \delta \Delta = \Img \Delta \delta 
\subset \Img \delta \cap \Ker \Delta$
and $\Img \Delta \delta = \Img  \delta \Delta
\subset \Img \Delta. \cap \Ker \delta$.
The following conditions are equivalent
\begin{itemize}
\item[(i)] The inclusions 
$i: (\Ker \Delta, \delta) \hookrightarrow (\cA, \delta)$
and $j: (\Ker \delta, \Delta) \hookrightarrow (\cA, \Delta)$
induce isomorphisms of homology.

\item[(ii)] We have equalities:
\begin{align*}
&&\Img \delta \Delta = \Img \Delta \delta = \Img \delta \cap \Ker \Delta, \\
&&\Img \delta \Delta = \Img \Delta \delta = \Img \Delta \cap \Ker \delta. 
\end{align*}
\end{itemize}
When the above conditions hold,
then both the cohomology groups in (i) are naturally isomorphic to 
$(\Ker \Delta \cap \Ker \delta)/\Img \delta \Delta$.
\end{lemma}

From now on, all the DGBV algebras will be assumed
to satisfy the conditions in Lemma \ref{lm:quasi}.
Denote by $\DGBV_q$ the subcategory of all DGBV algebras which satisfies
the conditions in Lemma \ref{lm:quasi}.

\begin{lemma} \label{lm:tech}
If $z \in \cA$ satisfies $\delta \Delta z = 0$,
then
$$z = h + \Delta u + \delta v$$
for some $h \in \Ker \delta \cap \Ker \Delta$, $u, v \in \cA$.
\end{lemma}

\begin{proof}
Since $\Delta z \in \Ker \delta \cap \Img \Delta = \Img \delta \Delta$,
so $\Delta z = \Delta \delta v$ for some $v \in \cA$.
Now $\Delta (z - \delta v) = 0$,
so $z - \delta v$ determines a class in $H(\cA, \Delta)$.
Since $H(\cA, \Delta) \cong (\Ker \delta \cap \Delta)/\Img \delta \Delta$,
there exist $u \in \cA$,
such that
$h = z - \delta v - \Delta u \in Ker \delta \cap \Ker \Delta$.
This completes the proof. 
\end{proof}

\begin{proposition} \label{prop:MC}
Assume that  $(\cA, \wedge, \delta, \Delta, [\cdot\bullet\cdot])$ is 
in $\DGBV_q$.
Then for any even cohomology class $[x] \in H(\cA, \delta)$,
there is a formal power series solution 
$x(t) = x_1 t + \cdots + x_n t^n + \cdots$ to (\ref{eqn:MC}),
such that $x_1 \in \Ker \Delta \cap \Ker \delta$ represents $[x]$,
and $x_n \in \Img \Delta$ for $n > 1$.
When $[x] = [1]$,
we can take $x_1 = 1$, $x_n = 0$, for $n > 1$.
\end{proposition}

\begin{proof}
The second statement is trivial.
The first can be proved by a standard argument modeled on the method of 
Tian \cite{Tia} and Todorov \cite{Tod}.
Rewrite (\ref{eqn:MC}) as a sequence of equations
\begin{eqnarray*}
&& \delta x_1 = 0, \\
&& \delta x_2 = - \frac{1}{2} [x_1 \bullet x_1], \\
&& \cdots \cdots \\
&& \delta x_n = - \frac{1}{2} \sum_{i + j = n} [x_i \bullet x_j], \\
&& \cdots \cdots
\end{eqnarray*}
By Lemma \ref{lm:quasi},
we can take $x_1 \in \Ker \delta \cap \Ker \Delta$ to represent $[x]$.
Suppose now we have found $x_1, \cdots, x_n$ with $x_j \in \Img \Delta$ for 
$1 < j \leq n$.
By (\ref{eqn:bracket}),
$\sum_{i + j = n + 1} [x_i \bullet x_j] \in \Img \Delta$.
Also we have the following standard calculation
\begin{eqnarray*}
&& \delta \sum_{i + j = n + 1} [x_i \bullet x_j] 
=  \sum_{i + j = n + 1} [\delta x_i \bullet x_j]
	- \sum_{i + j = n + 1} [x_i \bullet \delta x_j] \\
& = & -\frac{1}{2} \sum_{i + j = n + 1} \sum_{p + q = i}
	[[x_p \bullet x_q] \bullet x_j]
+ 	\frac{1}{2} \sum_{i + j = n + 1} \sum_{p + q = j}
	[x_i \bullet [x_p \bullet x_q]] \\
& = & -\frac{1}{2} \sum_{i + j = n + 1} \sum_{p + q = i}
	[[x_p \bullet x_q] \bullet x_j]
+ 	\frac{1}{2} \sum_{i + j = n + 1} \sum_{p + q = j}
	[[x_p \bullet x_q] \bullet x_i] = 0.
\end{eqnarray*}
Hence $\sum_{i + j = n + 1} [x_i \bullet x_j] 
\in \Ker \delta \cap \Img \Delta = \Img \delta \Delta$,
and so there exist $x_{n+1} \in \Img \Delta$,
such that 
$$\delta x_{n+1} 
= - \frac{1}{2} \sum_{i + j = n+1} [x_i \bullet x_j].$$
\end{proof}

\begin{proposition} \label{prop:quasi2}
The triple $(\cA[[t]], \delta_{x(t)}, \Delta)$ 
is in $\DGBV_q$.
\end{proposition}

\begin{proof}
We need to show 
\begin{eqnarray*}
&& \Img \delta_{x(t)} \cap \Ker \Delta 
	\subset \Img \delta \Delta = \Img \Delta \delta, \\
&& \Img \Delta \cap \Ker \delta_{x(t)} 
	\subset \Img \delta \Delta = \Img \Delta \delta.
\end{eqnarray*}
We prove the second inclusion first.
Assume that $y(t) = \Delta z(t)$, $\delta_{x(t)}y(t) = 0$.
We get a sequence of equations
\begin{eqnarray*}
&& \delta y_0 = 0, \\
&& \delta y_1 = - [x_1 \bullet y_0], \\
&& \cdots \cdots \\
&& \delta y_n = - \sum_{i + j = n} [x_i \bullet y_j], \\
&& \cdots \cdots
\end{eqnarray*}
where $y_n = \Delta z_n$.
Now $\delta \Delta z_0 = \delta y_0 = 0$,
by Lemma \ref{lm:tech},
$z_0 = h_0 + \Delta u_0 + \delta v_0$, 
where $h_0 \in \Ker \delta \cap \Ker \Delta$.
Hence 
$$y_0 = \Delta z_0 = \Delta \delta v_0 = \delta (- \Delta v_0).$$
And from
\begin{eqnarray*}
\delta\Delta z_1 = \delta y_1 = - [x_1 \bullet y_0]
= [x_1 \bullet \delta \Delta v_0] = \delta \Delta [x_1 \bullet v_0],
\end{eqnarray*}
we get $z_1 = h_1 + \Delta u_1 + \delta v_1 + [x_1 \bullet v_0] $.
Hence we have
\begin{eqnarray*}
y_1 = \Delta z_1 = \Delta \delta v_1 + \Delta [x_1 \bullet v_0] 
= - \delta \Delta v_1 - [x_1 \bullet \Delta v_0].
\end{eqnarray*}
By induction, we find that 
$$y_n = - \delta \Delta v_n - \sum_{i + j = n} [x_i \bullet v_j],$$
in other word, 
$y(t) = - \delta_{x(t)} \Delta v(t) \in \Img \delta_{x(t)} \Delta$,
where $v(t) = v_0 + v_1 t + \cdots$.

Now assume $y(t) = \delta_{x(t)} z(t)$ and $\Delta y(t) = 0$, 
where $y(t) = y_0 + y_1 t + \cdots$ and $z(t) = z_0 + z_1 t+ \cdots$
are elments of $\cA[[t]]$.
Equivalently, we have  a sequence of equations
\begin{eqnarray*}
&& y_0 = \delta z_0, \\
&& y_1 = \delta z_1 + [x_1 \bullet z_0], \\
&& \cdots \cdots \\
&& y_n = \delta z_n + \sum_{i + j = n} [x_i \bullet z_j], \\
&& \cdots \cdots 
\end{eqnarray*}
and $\Delta y_n = 0$.
Now  $\Delta \delta z_0 = \Delta y_0 = 0$,
$z_0 = h_0 + \Delta u_0 + \delta v_0$ 
for some $h_0 \in \Ker \delta \cap \Ker \Delta$.
Hence 
$$y_0 = \delta z_0 = \delta \Delta v_0.$$
From
\begin{eqnarray*}
&& \Delta \delta z_1 = \Delta (y_1 - [x_1 \bullet z_0]) 
= [x_1 \bullet \Delta z_0] 
= [x_1 \bullet \Delta \delta v_0]
= \Delta \delta [x_1 \bullet v_0],
\end{eqnarray*}
we get 
$z_1 = h_1 + \Delta u_1 + \delta v_1 + [x_1 \bullet v_0]$.
Hence
\begin{eqnarray*}
y_1 & = & \delta z_1 + [x_1 \bullet z_0]
= \delta \Delta u_1 + \delta [x_1 \bullet v_0] 
+ [x_1 \bullet (h_0 + \Delta u_0 + \delta v_0)] \\
& = & \delta \Delta u_1 + [x_1 \bullet h_0] + [x_1 \bullet \Delta u_0].
\end{eqnarray*}
By induction, 
we can show that
\begin{eqnarray*}
z_n = h_n + \Delta u_n + \delta v_n + \sum_{i+j= n} [x_i \cdot v_j],
\end{eqnarray*}
where each $h_n$ lies in $\Ker \delta \cap \Ker \Delta$.
Consequently,
\begin{eqnarray*}
y_n & = & \delta \Delta u_n + \sum_{i + j = n} [x_i \bullet h_j]
+ \sum_{i + j = n} [x_i \bullet \Delta u_j]
\end{eqnarray*}
Define $w_n$ as follows.
Set $w_0 = 0$.
For $n > 0$,
assume that $w_1, \cdots w_{n-1}$ have been defined
such that
$$\delta \Delta w_j = \sum_{q + r = j} [x_p \bullet (h_q - \Delta w_q)],$$
for $j < n$.
Then it is straightforward to see that
$\sum_{i + j = n} [x_i \bullet (h_j - \Delta w_j)]  
= \sum_{i + j = n} \Delta (x_i \wedge (h_j - \Delta w_j)) 
\in \Ker \delta \cap \Img \Delta = \Img \delta \Delta$,
so 
$$ \sum_{i + j = n} [x_i \bullet (h_j - \Delta w_j)]= \delta \Delta w_n$$ 
for some $w_n \in \cA$.
Finally, we have
\begin{eqnarray*}
y_n & = & \delta \Delta u_n + \sum_{i + j = n} [x_i \bullet (h_j - \Delta w_j]
+ \sum_{i + j = n} [x_i \bullet \Delta (u_j + w_j)] \\
& = & \delta \Delta u_n + \delta \Delta w_n 
	+ \sum_{i + j = n} [x_i \bullet \Delta (u_j + w_j)] \\
& = & \delta \Delta (u_n + w_n) 
	+ \sum_{i + j = n} [x_i \bullet \Delta (u_j + w_j)].
\end{eqnarray*}
\end{proof}

\begin{proposition} \label{prop:extension}
Any element $y_0 \in \Ker \delta \cap \Ker \Delta$ can be extended 
to a formal power series $y(t) = y_0 + y_1 t + \cdots y_n t^n + \cdots$,
such that 
$\delta_{x(t)} y(t) = 0$ and $y_n \in \Img \Delta$ for $n \geq 1$.
Furthermore, if $\bar{y}(t) = \bar{y}_0 + \bar{y}_1 t + \cdots$ satisfies
$\delta_{x(t)} \bar{y}(t) = 0$,
where $\bar{y}_0 \in \Ker \delta \cap \Ker \Delta$ represents the same class
as $y_0$ in $H(\cA, \delta)$,
$\bar{y}_n \in \Img \Delta$ for $n > 0$,
then there exists $z(t) \in (\Img \Delta)[[t]]$ such that
$\bar{y}(t) - y(t) = \delta_{x(t)} z(t)$.
\end{proposition}

\begin{proof}
When expanded into the formal power series in $t$,
we can rewrite $\delta_{x(t)} y(t) = 0$ as a sequence of equations
\begin{eqnarray*}
&& \delta y_0 = 0, \\
&& \delta y_1 = - [x_1 \bullet y_0], \\
&& \cdots\cdots \\
&& \delta y_n = - \sum_{i + j = n} [x_i \bullet y_j], \\
&& \cdots \cdots 
\end{eqnarray*}
To find $y_1$,
notice that 
\begin{eqnarray*}
&& [x_1 \bullet y_0] = \Delta (x_1 \wedge y_0) \in \Img \Delta, \\
&& \delta [x_1 \bullet y_0] 
	= [\delta x_1 \bullet y_0]
	+ [x_1 \bullet \delta y_0] = 0,
\end{eqnarray*}
i.e. $[x_1 \bullet y_0] 
\in \Ker \delta \cap \Img \Delta = \Img \delta \Delta$,
hence $y_1 \in \Img \Delta$ can be found.
Suppose now that we have found $y_1, \cdots, y_n \in \Img \Delta$,
we have
\begin{eqnarray*}
&& \sum_{i+j=n+1} [x_i \bullet y_j] 
	= \sum_{i+j=n+1}\Delta (x_i \wedge y_j) \in \Img \Delta, \\
&& \delta \sum_{i+j=n+1} [x_i \bullet y_j] 
	=  \sum_{i+j=n+1} [\delta x_i \bullet y_j]
	-  \sum_{i+j=n+1} [x_i \bullet \delta y_j] \\
& = & 	-\frac{1}{2}\sum_{i+j=n+1} \sum_{p+q = i}
	[[x_p \bullet x_q]\bullet y_j]
	+ \sum_{i+j=n+1} \sum_{q+ r = j} 
	[x_i \bullet [x_q \bullet y_r]] \\
& = & 	-\frac{1}{2} \sum_{p + q + r = n + 1}
([[x_p \bullet x_q]\bullet y_r]
- 2 [x_p \bullet [x_q \bullet y_r]]) = 0.
\end{eqnarray*}
I.e., $\sum_{i+j=n+1} [x_i \bullet y_j] \in \Ker \delta \Img \Delta 
= \Img \delta \Delta$,
hence we can find $y_{n+1}$.

By Lemma \ref{lm:quasi},
$H(\cA, \delta) \cong (Ker \Delta \cap \Ker \delta) / \Img \delta \Delta$,
since $\bar{y}_0$ and $y_0$ represents the same class in $H(\cA., \delta)$,
we have $\bar{y}_0 - y_0 = \delta z_0$ for some $z_0 \in \Img \Delta$. 
Now as above, we solve $\delta_{x(t)} z(t) = \bar{y}(t) - y(t)$ by induction:
first expand in power series to get a sequence of equations
\begin{eqnarray*}
&& \delta z_0 = \bar{y}_0 - y_0, \\
&& \delta z_1 = \bar{y}_1 - y_1 - [x_1 \bullet z_0], \\
&& \cdots \cdots \\
&& \delta z_n = \bar{y}_n - y_n - \sum_{i + j = n} [x_i \bullet z_j], \\
&& \cdots \cdots
\end{eqnarray*}
then inductively check the right hand side of each equation lies in
$\Ker \delta \cap \Img \Delta = \Img \delta \Delta$ as above,
hence one can find a solution $z_n$ in $\Img \Delta$.
\end{proof}

We now define a homomorphism 
$\phi_{x(t)}: H(\cA, \delta)[[t]] 
\to H(\cA[[t]], \delta_{x(t)})$
as follows: for any $[y] \in H(\cA, \delta)$,
represent it by an element $y_0 \in \Ker \delta \cap \Ker \Delta$.
Extend $y_0$ to $y(t)$ as in Proposition \ref{prop:extension},
then set $\phi_{x(t)}([y]) = [y(t)]$.
This is well-defined by Proposition \ref{prop:extension}.
Extend $\phi_{x(t)}$ to $H(\cA, \delta)[[t]]$
and still denote it by $\phi_{x(t)}$.
We have the following

\begin{proposition} \label{prop:isomorphism}
The homomorphism 
$\phi_{x(t)}:  H(\cA, \delta)[[t]] \to H(\cA[[t]], \delta_{x(t)})$
is an isomorphism of $\bk[[t]]$-modules.
\end{proposition}

\begin{proof}
We break the proof into two steps.

\noindent {\em Step 1. $\phi_{x(t)}$ is injective.}
By Lemma \ref{lm:quasi},
we can represent any class of $H(\cA, \delta)[[t]]$
by an element $y(t) = y^{(k)} t^k + y^{(k+1)} t^{k+1} + \cdots
\in (\Ker \delta \cap \Ker \Delta)[[t]]$. 
Without loss of generality,
we assume that $y^{(k)} \notin \Img \Delta \delta$
and hence $[y^{(k)}] \neq 0 \in H(\Ker \Delta, \delta)$.
In fact, 
if $y^{(k)} = \Delta \delta u_k$ for some $u_k \in \cA$,
we replace $y(t)$ by $y(t) - t^k \delta\Delta u_k$
and consider the $(k+1)$-th term.
Now if $\phi_{x(t)}(y(t)) = \delta_{x(t)} z(t)$ for
some $z(t) = z_0 + z_1 t + \cdots \in (\Ker \Delta)[[t]]$,
then by noticing that $\phi_{x(t)}(y(t))$ has leading term $y^{(k)}t^k$,
we get a sequence of equations
\begin{eqnarray*}
&& \delta z_0 = 0, \\
&& \delta z_1 + [x_1 \bullet z_0 ] = 0, \\
&& \cdots \cdots \\
&& \delta z_{k -1} + \sum_{i + j = k - 1} [x_i \bullet z_j] = 0, \\
&& \delta z_k + \sum_{i + j = k} [x_i \bullet z_j] = y^{(k)}.
\end{eqnarray*}
Now $\sum_{i + j = k} [x_i \bullet z_j] 
= \sum_{i + j = k} \Delta (x_i \wedge z_j) \in \Img \Delta$,
and
\begin{eqnarray*}
&& \delta \sum_{i + j = k} [x_i \bullet z_j] 
= \sum_{i + j = k} [\delta x_i \bullet z_j] 
	- \sum_{i + j = k} [x_i \bullet \delta  z_j] \\
& = & \sum_{i + j = k} [\delta x_i \bullet z_j] 
+ \sum_{i + j = k} [x_i \bullet \sum_{q + l = j}[x_q \bullet z_l]] \\
& = & \sum_{i + j = k} [\delta x_i \bullet z_j] 
+ \sum_{i + j = k} \frac{1}{2} \sum_{p+q = i} 
	[[x_p \bullet x_q] \bullet z_j] = 0.
\end{eqnarray*}
Therefore $\sum_{i + j = k} [x_i \bullet z_j] \in \Img \delta \Delta$,
and so $y^{(k)} \in \delta \Ker \Delta$,
a contradiction to the assumption 
that $[y^{(k)}] \neq 0 \in H(\Ker \Delta, \delta)$.

\noindent {\em Step 2. $\phi_{x(t)}$ is surjective.}
By Lemma \ref{prop:quasi2},
any element of $H(\cA[[t]], \delta_{x(t)})$ can be represented by
an element 
$y(t) = y^{(0)} + y^{(1)} t + \cdots\in \Ker \Delta \cap \Ker \delta_{x(t)}$.
Then we have $y^{(0)} \in \Ker \Delta \cap \Ker \delta$,
hence it can be extended to an element 
$\phi_{x(t)} (y^{(0)}) = y^{(0)} + y^{(0)}_1 t + \cdots \in
\Ker \Delta \cap \Ker \delta_{x(t)}$.
Consider now $y(t) - \phi_{x(t)} (y^{(0)})$,
it can be written as $t y'(t)$,
where $\tilde{y}(t) \in \Ker \Delta \cap \Ker \delta_{x(t)}$.
Hence by induction, $y(t)$ lies in the image of $\phi_{x(t)}$.
\end{proof}

\begin{theorem}
Assume that $(\cA, \wedge, \delta, \Delta, [\cdot \bullet \cdot])$ 
is a DGBV algebra in $\DGBV_q$ and
$x(t)$ as in Proposition \ref{prop:MC}.
Then there is a naturally defined formal deformation 
$$\wedge_{x(t)}:  H(\cA, \delta)[[t]] \otimes H(\cA, \delta)[[t]] \to
H(\cA, \delta)[[t]]$$ of $\wedge$
defined by
$a \wedge_{x(t)} b = \phi^{-1}_{x(t)}(\phi_{x(t)}(a) \wedge \phi_{x(t)}(b))$.
\end{theorem}

\section{Nice integrals and formal Frobenius manifold structures}

In the above discussion,
we have only considered one-parameter deformations such that the parameter
$t$ commutes with the elements of the relevant algebras.
In this section,
we need to make some generalizations.

First we will consider odd deformations. 
Let $\epsilon$ be an indeterminate,
denote by $\bk[[\epsilon]]$ 
the exterior algebra $\Lambda^*(\bk \oplus \bk \epsilon)$.
We regard $\bk[[\epsilon]]$ as the super formal power series with 
an indeterminate of odd degree.
If $(\cA, \wedge)$ is a $\bZ_2$-graded associative algebra over $\bk$,
let $\cA[[\epsilon]] = \cA \otimes_{\bk} \bk[[\epsilon]]$.
Any element of $\cA[[t]]$ can be written as $a_0 + a_1 \epsilon$,
and  we can extend $\wedge$ as follows:
$$(a_0 + a_1 \epsilon) \wedge (b_0 + b_1 \epsilon)
= a_0 \wedge b_0 + (a_0 \wedge b_1 + a_1 \sigma(b_0) \epsilon,$$
where $\sigma(a) = (-1)^{|a|} a$ for a homogeneous element $a \in \cA$. 
Here we have used the Koszul sign convention.
It is straightforward to verify that this multiplication is associative.
Furthermore, if $1$ is the unit for $(\cA, \wedge)$,
then so is it for $(\cA[[\epsilon]], \wedge)$;
and if $(\cA, \wedge)$ is graded commutative,
then so is $(\cA[[\epsilon]], \wedge)$.
Now if $(\cA, \wedge, \delta, \Delta, [\cdot \bullet\cdot])$ is a DGBV algebra,
we extend $\Delta$ and $[\cdot \bullet \cdot]$ to $\cA[[\epsilon]]$
by the Koszul sign convention.
We have
$\epsilon a = \sigma(a) \epsilon$,
$\epsilon \delta = - \delta \epsilon$, $\epsilon \Delta = -\Delta \epsilon$
and 
\begin{eqnarray*}
[a \epsilon \bullet b] = -[a \bullet \sigma(b)] \epsilon, &
[a \bullet b \epsilon] = [a \bullet b] \epsilon.
\end{eqnarray*}
Given a cohomology class of $H(\cA, \delta)$ of odd degree,
represent it by an element in $\Ker \delta \cap \Ker \Delta$.
Then $x(\epsilon) = x_1 \epsilon$ satisfies the Maurer-Cartan equation
$$\delta \omega + \frac{1}{2}[\omega \bullet \omega] = 0$$
over $\cA[[\epsilon]]$.
This is the analogue of Proposition \ref{prop:MC}.
The simpler versions of the proofs (without inductions) 
of Propositions \ref{prop:quasi2} - \ref{prop:isomorphism} 
prove the corresponding statements for odd cohomology classes 
of $H(\cA, \delta)$.
As a result, we obtain an odd deformation of $(H(\cA, \delta), \wedge)$. 

Secondly, 
we will also be interested in multi-parameter formal deformations 
which now turn to.
From now on,
assume that $H$ is a rank $n$ free $\bk$-module.
Let $\{e_{\alpha}: \alpha = 0, \cdots, n-1\}$ 
be a set of free homogeneous generators of $H$,
such that $e_0 = 1$.
Let $\{t^{\alpha} \}$ be the dual set of 
generators of $H^t = \Hom_{\bk}(H, \bk)$.
Denote by $\bk[[\bt]]$ the space of super power
series in $\{t^0, \cdots, t^{n-1}\}$,
and consider $\cA[\bt]] = \cA \otimes_{\bk} \bk[[\bt]]$.
Again, modifying the proofs,
the analogues of Propositions \ref{prop:quasi2} - \ref{prop:isomorphism} 
can be proved for $\cA[[\bt]]$.
We give the corresponding statements below,
and omit the proofs.
Also, see Manin \cite{Man2} for different proofs.

\begin{proposition} \label{prop:MC'}
There is a formal power series solution 
$\Gamma = \Gamma_1  + \cdots + \Gamma_n + \cdots$ to (\ref{eqn:MC}),
such that $\Gamma_1 = x_{\alpha} t^{\alpha}$,
where $x_{\alpha} \in \Ker \Delta \cap \Ker \delta$ represents $e_{\alpha}$,
and for $n > 1$, $\Gamma_n$ is a super  homogeneous polynomial of order $n$ in 
$\{t^{\alpha}\}$ with coefficients in $\Img \Delta$ .
Furthermore, $t^0$ only appears in $\Gamma_1$.
(Such a solutions is called a {\em universal normalized solution}.)
\end{proposition}

\begin{proposition} \label{prop:quasi2'}
Let $\Gamma$ be a universal normalized solution.
Then the deformed triple $(\cA[[\bt]], \delta_{\Gamma}, \Delta)$ 
is in $\DGBV_q$.
\end{proposition}

\begin{proposition} \label{prop:extension'}
Let  $\Gamma$ be a universal normalized solution.  
Then any element $y_0 \in \Ker \delta \cap \Ker \Delta$ can be extended 
to $y(\bt) \in \cA[[\bt]]$ with leading term $y_0$,
such that 
$\delta_{\Gamma} y(t) = 0$ and all the higher order terms lies in $\Img \Delta$.
Furthermore, if $\bar{y}(\bt) \in \cA[[\bt]]$ has the same properties and 
its leading term $\bar{y}_0 \in \Ker \delta \cap \Ker \Delta$ 
represents the same class as $y_0$ in $H(\cA, \delta)$,
then there exists $z(\bt) \in (\Img \Delta)[[\bt]]$ such that
$\bar{y}(t) - y(t) = \delta_{\Gamma} z(t)$.
\end{proposition}

Similar to the definition of $\phi_{x(t)}$,
we define a homomorphism 
$\phi_{\Gamma}: H(\cA, \delta)[[\bt]] 
\to H(\cA[[t]], \delta_{\Gamma})$.

\begin{proposition} \label{prop:isomorphism'}
The homomorphism 
$\phi_{\Gamma}:  H(\cA, \delta)[[\bt]] \to H(\cA[[\bt]], \delta_{\Gamma})$
is an isomorphism of $\bk[[\bt]]$-modules.
\end{proposition}

We can now define an $n$-parameter super formal deformation
$$\wedge_{\Gamma}:  H(\cA, \delta)[[\bt]] \otimes H(\cA, \delta)[[\bt]] \to
H(\cA, \delta)[[\bt]]$$
by
$a \wedge_{\Gamma} b 
= \phi^{-1}_{\Gamma}(\phi_{\Gamma}(a) \wedge \phi_{\Gamma}(b))$.
Recall the following

\begin{definition}
An {\em integral} of a DGBV algebra 
$(\cA, \wedge, \delta, \Delta, [\cdot \bullet \cdot])$ 
is a $\bk$-linear $\int: \cA \to \bk$,
such that 
\begin{eqnarray*}
&& \int \delta a \wedge b = (-1)^{|a| + 1} \int a \wedge \delta b, \\
&& \int \Delta a \wedge b = (-1)^{|a|} \int a \wedge \Delta b,
\end{eqnarray*}
for any homogeneous $a, b \in \cA$.
Obviously, 
an integral on $\cA$ induces a well-defined supersymmetric bilinear form,
$g: H \otimes H \to \bk$ by
$$g([a], [b]) = \int a \wedge b.$$
An integral is called {\em nice}
if $g$ induces an isomorphism 
$H \to H^t$ by $[a] \in H \mapsto g([a], \cdot) \in H^t$.
\end{definition}

In the context of \S \ref{sec:deform},
we extend $\int$ to $\cA[[t]]$.
Then it is straightforward to see that $\int$ is an integral of 
$(\cA[[t]], \wedge, \delta_{x(t)}, \Delta, [\cdot\bullet\cdot])$
(Manin \cite{Man2}, Proposition 5.5.1).
If for $j = 1, 2$, $y_j \in \Ker \delta \cap \Ker \Delta$ such that
$\delta_{x(t)} (y_j + \Delta z_j(t)) = 0$ for some $z_j(t) \in \cA[[t]]$,
then we have
\begin{eqnarray*}
&& \int (y_1 + \Delta z_1(t)) \wedge (y_2 + \Delta z_2(t)) \\
& = & \int y_1 \wedge y_2 + \int y_1 \wedge \Delta z_2(t) + 
\int \Delta z_1(t) \wedge (y_2  + \Delta z_2(t)) \\
& = & \int y_1 \wedge y_2 +  \int \Delta \sigma(y_1) \wedge z_2(t)
+ \int \sigma(z_2(t)) \wedge \Delta (y_2  + \Delta z_2(t)) \\
& = & \int y_1 \wedge y_2.
\end{eqnarray*}
So if the integral $\int$ is nice for $\cA$,
so is it for $\cA[[t]]$.
Given a nice integral $\int$ on a DGBV algebra $\cA$,
since we have
$$g([a] \wedge [b], [c]) = \int ([a] \wedge [b]) \wedge [c] =
\int [a] \wedge ([b] \wedge [c]) = g([a], [b] \wedge [c]),$$
for $[a], [b], [c] \in H$,
$(H, \wedge, g)$ is a {\em (graded) Frobenius algebra}.
Therefore, given a nice integral of $\cA$, 
$(H[[t]], \wedge_{x(t)}, g)$ is a formal deformation of the Frobenius algebra
$(H, \wedge, g)$. 
Similarly, given a universal normalized solution $\Gamma$,
the above discussion can be carried out for 
$(\cA[[\bt]], \delta_{\Gamma}, \wedge, \delta, [\cdot\bullet\cdot])$,
a nice integral on $\cA$ then gives a multi-parameter formal deformation
$(H[[t]], \wedge_{\Gamma}, g)$ of $(H, \wedge, g)$.
We now state the following result in 
Barannikov-Kontsevich \cite{Bar-Kon} and Manin \cite{Man2}:

\begin{proposition} \label{prop:Frobenius}
Let $(\cA, \wedge, \delta, \Delta, [\cdot\bullet\cdot])$
be a DGBV algebra over $\bk$ in $\DGBV_q$.
Assume that $\bk$-module $H = H(\cA, \delta)$ has free homogeneous generators
$\{e_{\alpha}: \alpha = 0, \cdots, n -1\}$ with $e_0 = 1$,
and $\{ t^{\alpha}\}$ the dual generators of the dual module $H^t$.
Also assume that there is a nice integral $\int$ on $\cA$.
Then given a universal normalized solution $\Gamma = \Gamma_1 + \Delta B$
with $\Gamma_1 = e_{\alpha} t^{\alpha}$,
the super formal power series
\begin{eqnarray} \label{eqn:potential}
\Phi = \int (\frac{1}{6} \Gamma^3 - \frac{1}{2}\delta B \wedge \Delta B)
= \int (\frac{1}{6} \Gamma^3 - \frac{1}{4} (\Gamma - \Gamma_1) \wedge \Gamma^2)
\end{eqnarray}
satisfies 
\begin{eqnarray*}
\frac{\partial^3 \Phi}{\partial t^0 \partial t^{\alpha} t^{\beta}}
= g_{\alpha\beta},
\end{eqnarray*}
and the WDVV (Witten-Dijkgraaf-Verlinde-Verlinde) equations
\begin{eqnarray*}
\frac{\partial^3 \Phi}{\partial t^{\alpha} \partial t^{\beta} \partial t^{\mu}}
g^{\mu \nu} 
\frac{\partial^3 \Phi}{\partial t^{\nu} \partial t^{\gamma} \partial t^{\delta}}
= (-1)^{|e_{\alpha}|(|e_{\beta}| + |e_{\gamma}|)}
\frac{\partial^3 \Phi}{ \partial t^{\beta} \partial t^{\gamma} \partial t^{\mu}}
g^{\mu \nu} 
\frac{\partial^3 \Phi}{\partial t^{\nu} \partial t^{\alpha} \partial t^{\delta}},
\end{eqnarray*}
where $g_{\alpha\beta} = g(e_{\alpha}, e_{\beta}) 
= \int e_{\alpha} \wedge e_{\beta}$ and $(g^{\alpha\beta})$ is
the inverse matrix of $(g_{\alpha\beta})$.
\end{proposition}

\section{Gauge invariance and functorial property}

In this section, 
we will discuss the gauge invariance of 
the constructions in previous sections.
As mentioned in the introduction, 
special cases have been treated in the 
Kodaira-Spencer theory gravity \cite{Ber-Cec-Oog-Vaf} 
and the K\"{a}hler gravity \cite{Ber-Sad}.
Barannikov and Kontsevich \cite{Bar-Kon}
remarked that all universal normalized solutions are gauge equivalent 
but offered no proof.

Let $(L = \oplus_{n \in \bZ} L_n, [\cdot, \cdot], \delta)$ 
be a differential graded Lie algebra (DGLA).
Denote by $L^e$ and $L^o$ the subspaces of elements with even and odd 
degrees respectively.
For any element $A \in L$, denote by $\ad_A$ 
the automorphism $[A, \cdot]: L \to L$. 
Let $R = \bk[t]$ or $\bk[[\bt]]$,
and $\fm$ the maximal ideal of $R$.
Consider the group
\begin{eqnarray*}
\cG = G(R) = \exp (L \otimes \fm)^e,
\end{eqnarray*}
with the multiplication $e^A e^B = e^C$  defined
by the Campbell-Baker-Hausdorff formula: 
$C = \sum_{n \geq 1} C_n$, 
where $C_n = \frac{1}{n} \sum_{p + q = n} (C_{p, q}' + C_{p, q}'')$,
and
\begin{eqnarray*}
&& C_{p, q}' = \sum_{\substack{p_1 + \cdots p_m = p \\
q_1 + \cdots + q_{m-1} = q - 1 \\
p_i + q_i \geq 1 \\
p_m \geq 1}} 
\frac{(-1)^{m+1}}{m}
\frac{\ad_A^{p_1}\ad_B^{q_1} \cdots \ad_A^{p_m}B}
{(p_1)!(q_1)!\cdots (p_m)!} \\
&& C_{p, q}'' = \sum_{\substack{p_1 + \cdots p_{m-1} = p -1 \\
q_1 + \cdots + q_m = q \\
p_i + q_i \geq 1}} 
\frac{(-1)^{m+1}}{m}
\frac{\ad_A^{p_1}\ad_B^{q_1} \cdots \ad_B^{q_{m-1}}A}
{(p_1)!(q_1)!\cdots (p_m)!}
\end{eqnarray*}
These are the explicit formulas in Dynkin's form (Serre \cite{Ser}, p. 29). 
(This is similar to a construction in 
Goldman and Millson \cite{Gol-Mil1, Gol-Mil2}
where they use Artin local $\bk$-algebras.) 
There is a natural action of $\cG$ on $L \otimes R$ by
\begin{eqnarray} \label{eqn:action}
e^A \cdot \alpha = e^{\ad_A} \alpha.
\end{eqnarray}
It is clear that $e^{\ad_A} e^{-\ad_A} = 1$.

\begin{lemma} \label{lm:commutator}
Let $A \in (L \otimes \fm)^e$,
$B \in L \otimes R$,
then for $n \geq 1$, we have
\begin{eqnarray*}
&& \ad_B \ad_A^n 
= \sum_{j = 0}^n \left( \begin{array}{c} n \\ j \end{array} \right) 
\ad_A^j \ad_{(-\ad_A)^{n-j}B}.
\end{eqnarray*}
Furthermore, we have
\begin{eqnarray*}
e^{\ad_A}\ad_B e^{-\ad_A} = \ad_{e^{\ad_A} B}.
\end{eqnarray*}
\end{lemma}

\begin{proof}
The first equality can be proved elementarily by induction.
For the second equality, we have
\begin{eqnarray*}
&& \ad_B e^{\ad_A} = \sum_{n \geq 0} \frac{1}{n!} \ad_B \ad_A^n 
 	=  \sum_{n \geq 0} \frac{1}{n!} 
\sum_{j = 0}^n \left( \begin{array}{c} n \\ j \end{array} \right) 
\ad_A^j \ad_{(-\ad_A)^{n-j}B} \\
& = & \sum_{n \geq 0} 
\sum_{j = 0}^n \frac{1}{j!} \ad_A^j \frac{1}{(n-j)!} \ad_{(-\ad_A)^{n-j}B} 
= \sum_{j, k \geq 0} \frac{1}{j!} \ad_A^j \frac{1}{k!} \ad_{(-\ad_A)^k B} \\
& = & e^{\ad_A} \ad_{e^{-\ad_A}B}. 
\end{eqnarray*}
\end{proof}

\begin{lemma} \label{lm:commutator2}
For $A \in (L \otimes \fm)^e$, we have
$$e^{\ad_A} d e^{- \ad_A} 
= d + \sum_{q \geq 0} \frac{1}{(q+1)!} \ad_{\ad_A^q \delta A}
=d + \ad_{\frac{1 - e^{\ad_A}}{\ad_A} \delta A}.$$
\end{lemma}

\begin{proof}
It is easy to see that
\begin{eqnarray*}
d \ad_A = \ad_A d + \ad_{\delta A}. 
\end{eqnarray*}
Then by induction,
it is easy to show that
\begin{eqnarray*}
&& d \ad_A^n 
= \ad_A^n d + \sum_{p=0}^{n-1} \ad_A^{n - p - 1} \ad_{\delta A} \ad_A^p
\end{eqnarray*}
Now $\delta A \in L^o$,
so we can use Lemma \ref{lm:commutator} to handle
$\ad_{\delta A} \ad_A^p$ as follows:
\begin{eqnarray*}
&& d \ad_A^n 
= \ad_A^n d + \sum_{p=0}^{n-1} \ad_A^{n - p - 1} \ad_{\delta A} \ad_A^p \\
& = & \ad_A^n d + \ad_A^{n - 1} \ad_{\delta A} 
+ \sum_{p=1}^{n-1} \ad_A^{n - p - 1} \ad_{\delta A} \ad_A^p \\
& = & \ad_A^n d + \ad_A^{n - 1} \ad_{\delta A} 
+ \sum_{p=1}^{n-1} \ad_A^{n - p - 1} 
\sum_{j = 0}^p \left( \begin{array}{c} p \\ j \end{array} \right) 
\ad_A^j \ad_{(-\ad_A)^{p-j}\delta A} \\
& = & \ad_A^n d + \ad_A^{n - 1} \ad_{\delta A} 
+ \sum_{p=1}^{n-1} 
\sum_{j = 0}^p \left( \begin{array}{c} p \\ j \end{array} \right) 
\ad_A^{n - p - 1 + j} \ad_{(-\ad_A)^{p-j}\delta A} \\
& = & \ad_A^n d + n \ad_A^{n - 1} \ad_{\delta A} 
+ \sum_{p=1}^{n-1} 
\sum_{j = 0}^{p-1} \left( \begin{array}{c} p \\ j \end{array} \right) 
\ad_A^{n - p - 1 + j} \ad_{(-\ad_A)^{p-j}\delta A} \\
& = & \ad_A^n d + n \ad_A^{n - 1} \ad_{\delta A} 
+ \sum_{q=1}^{n-1} 
\sum_{j = 0}^{n-1-q} \left( \begin{array}{c} q+j \\ j \end{array} \right) 
\ad_A^{n - 1 - q} \ad_{(-\ad_A)^q \delta A} \\
& = & \ad_A^n d + n \ad_A^{n - 1} \ad_{\delta A} 
+ \sum_{q=1}^{n-1} \left( \begin{array}{c} n \\ n-1-q \end{array} \right) 
\ad_A^{n - 1 - q} \ad_{(-\ad_A)^q \delta A} \\
& = & \ad_A^n d  
+ \sum_{q=0}^{n-1} \left( \begin{array}{c} n \\ n-1-q \end{array} \right) 
\ad_A^{n - 1 - q} \ad_{(-\ad_A)^q \delta A}
\end{eqnarray*}
Hence we have
\begin{eqnarray*}
&& d e^{\ad_A} = \sum_{n \geq 0} \frac{1}{n!} d \ad_A^n \\
& = & \sum_{n \geq 0} \frac{1}{n!} \ad_A^n d  
+ \sum_{n \geq 0} \frac{1}{n!} 
\sum_{q=0}^{n-1} \left( \begin{array}{c} n \\ n-1-q \end{array} \right) 
\ad_A^{n - 1 - q} \ad_{(-\ad_A)^q \delta A} \\
& = & e^{\ad_A} d + \sum_{n \geq 0} \sum_{q=0}^{n-1}  
\frac{1}{(n-1-q)!}\ad_A^{n - 1 - q} 
\frac{1}{(q+1)!} \ad_{(-\ad_A)^q \delta A} \\ 
& = & e^{\ad_A} 
(d + \sum_{q \geq 0} \frac{1}{(q+1)!} \ad_{(-\ad_A)^q \delta A}).
\end{eqnarray*}
Replacing $A$ by $-A$ then completes the proof.
\end{proof}

Given an element $\omega \in (L \otimes \fm)^o$,
set $\delta_\omega = \delta + \ad_{\omega}$.
From Lemma \ref{lm:commutator} and Lemma \ref{lm:commutator2},
we see that for any $A \in (L\otimes \fm)^e$,
we have
\begin{eqnarray*} 
e^A \delta_{\omega} e^{-A} = \delta_{e^A \cdot \omega},
\end{eqnarray*}
where
\begin{eqnarray*}
e^A \cdot \omega = e^{\ad_A} \omega 
+ \frac{1 - e^{\ad_A}}{\ad_A} \delta A
\end{eqnarray*}
is called the gauge transformation of $\omega$.
This formula is well-known.
See e.g. Goldman-Millson \cite{Gol-Mil1, Gol-Mil2}.
Let $X = \{ \omega \in (L \otimes \fm)^o: 
\delta\omega + \frac{1}{2}[\omega, \omega] =0 \}$,
then the $\cG$-action preserves $X$.

Given any DGBV algebra $(\cA, \wedge, \delta, \Delta, [\cdot \bullet \cdot])$,
denote by $\cA[-1]$ the graded vector space with $(\cA[-1])_p = \cA_{p+1}$.
Then $(\cA[-1], \delta, [\cdot \bullet \cdot])$ is a DGLA.
When $A = \bk[[t]]$,
we get the group $\cG[[t]] = \{\exp y(t): y(t) \in t\cA^o[[t]] \}$
and its  natural action on $\cA[[t]]$.

\begin{proposition} \label{prop:gauge1}
Assume that $x(t) = x_1 t + \cdots + x_n t^n + \cdots$ and 
$\bar{x}(t) = \bar{x}_1 t + \cdots + \bar{x}_n t^n + \cdots$ are two solutions 
of the Maurer-Cartan equation
$$\delta \omega + \frac{1}{2} [\omega \bullet \omega] = 0,$$
such that
$x_1, \bar{x}_1 \in \Ker \delta \cap \Ker \Delta$,
and $x_n, \bar{x}_n \in \Img \Delta$.
If $[x_1] = [\bar{x}_1]$,
then there exists $w(t) = w_1 t + \cdots + w_n t^n + \cdots \in \cA[[t]]$ 
such that
$\exp (\Delta w(t)) \cdot x(t) = \bar{x}(t)$.
\end{proposition}

\begin{proof}
It suffices to prove the following:
if $x(t) = \bar{x}(t) \pmod {t^n}$, $n \geq 1$,
then there exist $z_n \in \cA$,
such that
$\exp (t^n \Delta z_n) \cdot x(t) = \bar{x}(t) \pmod {t^{n+1}}$.

For $n = 1$,
we clear have $x(t) = \bar{x}(t) \pmod t$.
Since $[\bar{x}_1] = [x_1]$,
there exists $z_1$ such that 
$\bar{x}_1 = x_1 - \delta \Delta z_1$.
Now modulo $t^2$, we have
\begin{eqnarray*}
\exp (t \Delta z_1) \cdot x(t) 
& = & x(t) + \sum_{j > 0} \frac{1}{j!} \ad_{t \Delta z_1}^j x(t) 
- \delta \Delta z_1 - \sum_{j > 1} \frac{1}{j!} \ad_{t \Delta z_1}^{j-1} x(t) \\
& = & x_1 t - \delta \Delta z_1 = \bar{x}_1 t =\bar{x}(t).
\end{eqnarray*}
For $n > 1$,
we have $\bar{x}_n - x_n \in \Img \Delta$,
and
\begin{eqnarray*}
&& \delta (\bar{x}_n - x_n) = \delta \bar{x}_n - \delta x_n \\
& = & \sum_{i + j = n} [\bar{x}_i \bullet \bar{x}_j] 
	- \sum_{i+j =n} [x_i \bullet x_j] \\
& = & \sum_{i+j =n} [x_i \bullet x_j] - \sum_{i+j =n} [x_i \bullet x_j] = 0,
\end{eqnarray*}
hence there exists $z_n$ such that
$\bar{x}_n - x_n = - \delta \Delta z_n$.
Now modulo $t^{n+1}$, we have
\begin{eqnarray*}
&& \exp (t^n \Delta z_n) \cdot x(t) \\
& = & x(t) + \sum_{j > 0} \frac{1}{j!} \ad_{t^n \Delta z_n}^j x(t) 
- t^n \delta \Delta z_n 
- \sum_{j > 1} \frac{1}{j!} \ad_{t^n \Delta z_n}^{j-1} x(t) \\
& = & x_1 t + \cdots + x_{n-1}t^{n-1} + (x_n - \delta \Delta z_n) t^n 
= \bar{x}_1 t + \cdots \bar{x}_n t^n = \bar{x}(t).
\end{eqnarray*}
\end{proof}

As an easy corollary,
we have

\begin{proposition}
The equivalence class of $\wedge_{x(t)}$ depends only on the class $[x]$.
\end{proposition}

Consequently,
we have the following

\begin{theorem}
If $f: (\cA_1, \wedge_1, \delta_1, \Delta_1, [\cdot\bullet\cdot]_{\Delta_1}) 
\to
(\cA_2, \wedge_2, \delta_2, \Delta_2, [\cdot\bullet\cdot]_{\Delta_2})$
is a homomorphism of two DGBV algebras in $\DGBV_q$,
then we have
$$f_*(\alpha \wedge_x \beta) 
= f_*(\alpha) \wedge_{f_*(x)} f_*(\beta),$$
where $\alpha, \beta, x \in H(\cA_1, \delta_1)$,
and $f_*: H(\cA_1, \delta_1) \to H(\cA_2, \delta_2)$ 
is the homomorphism on cohomology induced by $f$. 
\end{theorem}
 
\begin{proof}
Given any $x \in H(\cA_1, \delta_1)$,
represent it by an element $x_1 \in \Ker \delta_1 \cap \Ker \Delta_1$
and extend it to a power series $x(t) = x_1 t+ \cdots + x_n t^n + \cdots$,
such that $x_n \in \Img \Delta_1$ and 
$\delta x(t) + \frac{1}{2}[x(t)\bullet x(t)]_{\Delta_1} = 0$.
Then $f(x(t)) = f(x_1) t + \cdots + f(x_n) t^n + \cdots$
satisfies   $f(x_1) \in \Ker \delta_2 \cap \Ker \Delta_2$,
$f(x_n) \in \Img \Delta_2$ and 
$\delta f(x(t)) + \frac{1}{2}[f(x(t))\bullet f(x(t))]_{\Delta_2} = 0$.
\end{proof}

\begin{proposition} \label{prop:gauge2}
Assume $(\cA, \wedge, \delta, \Delta, [\cdot\bullet\cdot])$ is a DGBV algebra
which satisfies the conditions in Proposition \ref{prop:Frobenius}.
Then given any two universal normalized solutions 
$\Gamma$ and $\overline{\Gamma}$,
there exists an odd element $A \in (\Img \Delta)[[\bt]]$
whose zero order term vanishes, such that
$e^A \cdot \Gamma = \overline{\Gamma}$.
Furthermore,
the potential function $\Phi$ is gauge invariant.
\end{proposition}

\begin{proof}
The proof of the first statement is an easy modification of the proof of
Proposition \ref{prop:gauge1}.
To prove the second statement,
we first linearize the gauge transformation:
\begin{eqnarray*}
\left.\frac{d}{d \lambda}\right|_{\lambda =0}e^{\lambda A} \cdot \Gamma 
= [A \bullet \Gamma] - \delta A.
\end{eqnarray*}
Then we have
\begin{eqnarray*}
&&  \left.\frac{d}{d \lambda}\right|_{\lambda =0}
\Phi(e^{\lambda A} \cdot \Gamma) \\
& = & \left.\frac{d}{d \lambda}\right|_{\lambda =0}
\int \left(\frac{1}{6} (e^{\lambda A} \cdot \Gamma)^3 
- \frac{1}{4} (e^{\lambda A} \cdot \Gamma - (e^{\lambda A} \cdot \Gamma)_1) 
\wedge (e^{\lambda A} \cdot \Gamma)^2 \right)\\
& = & \int ( \frac{1}{2} \Gamma^2 \wedge ([A \bullet \Gamma] - \delta A)
 - \frac{1}{4}([A \bullet \Gamma] - \delta A - ([A \bullet \Gamma] - \delta A)_1) 
 	\wedge \Gamma^2 \\
&& - \frac{1}{2}(\Gamma - \Gamma_1) \wedge \Gamma 
\wedge ([A \bullet \Gamma] - \delta A)) \\
& = & \int ( -\frac{1}{4} \Gamma^2 \wedge ([A \bullet \Gamma] - \delta A)
+ \frac{1}{4} \Gamma^2 \wedge ([A \bullet \Gamma] - \delta A)_1 
+ \frac{1}{2} \Gamma_1 \wedge  \Gamma \wedge ([A \bullet \Gamma] - \delta A)),
\end{eqnarray*}
where the subscript means the the first order term.
Since $A \in (\Img \Delta)[[\bt]]$ is odd,
we have
\begin{eqnarray*}
&& \int \Gamma^2 \wedge [A \bullet \Gamma] 
= \frac{1}{3} \int [A \bullet \Gamma^3]
= - \frac{1}{3} \int (\Delta (A \wedge \Gamma^3) - \Delta A \wedge \Gamma^3
+ A \wedge \Delta \Gamma^3)  \\
& = & \frac{2}{3} \int \Delta A \wedge \Gamma^3 = 0.
\end{eqnarray*}
Notice that
\begin{eqnarray*}
\int \Gamma^2 \wedge \delta A =  - \int \delta \Gamma^2 \wedge A
= - 2 \int \Gamma \wedge \delta \Gamma \wedge A
= \int \Gamma \wedge [\Gamma \bullet \Gamma] \wedge A
\end{eqnarray*}
Now
\begin{eqnarray*}
&& \int \Gamma \wedge [\Gamma \bullet \Gamma] \wedge A
= \frac{1}{2} \int [\Gamma \bullet \Gamma^2] \wedge A \\
& = & \frac{1}{2} \int (\Delta \Gamma^3 - \Delta \Gamma \wedge \Gamma^2 - 
\Gamma \wedge \Delta \Gamma^2) \wedge A \\
& = & \frac{1}{2} \int (\Gamma^3 \wedge \Delta A 
- \Gamma \wedge [\Gamma \bullet \Gamma] \wedge A)
= - \frac{1}{2} \int \Gamma \wedge [\Gamma \bullet \Gamma] \wedge A,
\end{eqnarray*}
hence $\int \Gamma \wedge [\Gamma \bullet \Gamma] \wedge A = 0$,
and so 
$$\int \Gamma^2 \wedge \delta A = 0.$$
Replacing $A$ by $A_1$,
we get
$$\int \Gamma^2 \wedge \delta A_1 = 0.$$
Similarly, we have
\begin{eqnarray*}
&& \int \Gamma_1 \wedge \Gamma \wedge [A \bullet \Gamma] 
= \frac{1}{2} \int \Gamma_1 \wedge [A \bullet \Gamma^2]  \\
& = & - \frac{1}{2} \int \Gamma_1 \wedge (\Delta (A \wedge \Gamma^2)
- \Delta A \wedge \Gamma^2 + A \wedge \Delta \Gamma^2) \\
& = & - \frac{1}{2} \int \Delta \Gamma_1 \wedge (A \wedge \Gamma^2)
- \frac{1}{2} \int  \Gamma_1 \wedge A \wedge [\Gamma \bullet \Gamma] \\
& = & \int \Gamma_1 \wedge A \wedge \delta \Gamma 
= \int \delta (\Gamma_1 \wedge A) \wedge \Gamma \\
& = & \int (\delta \Gamma_1 \wedge A \wedge \Gamma 
	+ \Gamma_1 \wedge \delta A \wedge \Gamma) 
= \int \Gamma_1 \wedge \Gamma \wedge \delta A,
\end{eqnarray*}
and so
$$\int \Gamma_1 \wedge \Gamma \wedge ([A \bullet \Gamma] - \delta \Gamma) = 0.$$
Therefore, we have
\begin{eqnarray*}
\left.\frac{d}{d \lambda}\right|_{\lambda =0}
\Phi(e^{\lambda A} \cdot \Gamma)
= \int \Gamma^2 \wedge [A \bullet \Gamma]_1 
= 0.
\end{eqnarray*}
Here we use the observation that
since both $A$ and $\Gamma$ has no zeroth order term,
$[A \bullet \Gamma]$ has no first order term.
\end{proof}

\begin{definition} 
Denote by $\DGBV_{qi}$ the subcategory of $\DGBV_q$ consists of 
DGBV elements in $\DGBV$ which admits nice integrals.
A {\em quasi-isomorphism} between two elements in $\DGBV_q$ ($\DGBV_{qi}$)
is a morphism in $\DGBV_q$ ($\DGBV_{qi}$) which
induces isomorphism on cohomology groups.
\end{definition}

\begin{theorem}
Formal Frobenius manifolds obtained from quasi-isomorphic DGBV algebras 
in $\DGBV_{qi}$ by Proposition 2.5 can be identified with each other. 
\end{theorem}

\begin{proof}
Just observe that a universal normalized solution is mapped to a universal
normalized solution under quasi-isomorphism.
\end{proof}


\begin{thebibliography}{99}


\bibitem{Bar-Kon} S. Barannikov, M. Kontsevich,
{\em Frobenius Manifolds and Formality of 
Lie Algebras of Polyvector Fields},
preprint, alg-geom/9710032.


\bibitem{Ber-Cec-Oog-Vaf}
M. Bershadsky, S. Cecotti, H. Ooguri, C. Vafa, 
{\em Kodaira-Spencer theory of gravity and exact results for quantum string
amplitudes}
{\bf Comm. Math. Phys. 165} (1994), no. 2, 311--427.

\bibitem{Ber-Sad}
M. Bershadsky, V. Sadov, 
{\em Theory of K\"{a}hler gravity}, 
{\bf Internat. J. Modern Phys. A 11} (1996), no. 26, 4689--4730.

\bibitem{Cao-Zho1}
H-D. Cao, J. Zhou,
{\em Frobenius Manifold Structure on Dolbeault Cohomology and Mirror Symmetry},
to appear in {\bf Comm. Anal. Geom.}, math.DG/9805094.

\bibitem{Cao-Zho2}
H-D. Cao, J. Zhou,
{\em Identification of Two Frobenius Manifolds},
to appear in {\bf Math. Res. Lett.}, math.DG/9805095.

\bibitem{Cao-Zho3}
H-D. Cao, J. Zhou,
{\em DGBV algebras and Frobenius Manifolds from hyperk\"{a}hler manifolds},
preprint.

\bibitem{Cao-Zho4}
H-D. Cao, J. Zhou,
{\em Formal Frobenius manifold structure on equivariant cohomology },
preprint.

\bibitem{Dub1} B. Dubrovin,
{\em Integrable systems in topological field theory},
{\bf  Nuclear Phys. B 379} (1992), no. 3,
627--689. 

\bibitem{Dub2} B. Dubrovin,
{\em Geometry of $2$D topological field theories},
in {\bf Integrable systems and quantum groups}
(Montecatini Terme, 1993), 120--348, 
Lecture Notes in Math., 1620, Springer, Berlin, 1996. 

\bibitem{Gol-Mil1}
W.M. Goldman, J.J. Millson, 
{\em The deformation theory of representations of fundamental groups 
of compact K\"{a}hler manifolds},
{\bf  Inst. Hautes \'{E}tudes Sci. Publ. Math. No. 67} (1988), 43--96. 

\bibitem{Gol-Mil2}
W.M. Goldman, J.J. Millson, 
{\em The homotopy invariance of the Kuranishi space},
{\bf Illinois J. Math. 34} (1990), no. 2, 337--367.

\bibitem{Ler-Vaf-War}
W. Lerche, C. Vafa, N.P. Warner, 
{\em Chiral rings in $N=2$ superconformal theories},
{\bf Nuclear Phys. B 324} (1989), no. 2, 427--474.

\bibitem{Lia-Liu-Yau}
B.H. Lian, K. Liu, S.-T. Yau, 
{\em Mirror principle. I.},
{\bf Asian J. Math. 1} (1997), no. 4, 729--763. 
  

\bibitem{Mor1}
D.R. Morrison,
{\em Mirror symmetry and rational curves on quintic threefolds: 
a guide for mathematicians},
{\bf J. Amer. Math. Soc. 6} (1993), no. 1, 223--247.

\bibitem{Mor2}
D.R. Morrison,
{\em Picard-Fuchs equations and mirror maps for hypersurfaces.},
in {\bf Essays on mirror manifolds}, 241--264, Internat. Press,
Hong Kong, 1992.


\bibitem{Man1} Yu. Manin, 
{\em Frobenius manifolds, quantum cohomology,
and moduli spaces (Chapters I, II, III)}.
Preprint MPI 96-113, 1996.


\bibitem{Man2} Y. Manin, 
{\em Three constructions of Frobenius manifolds: a comparative study},
preprint, math.QA/9801006.

\bibitem{Ser}
J.-P. Serre, 
{\bf Lie algebras and Lie groups}. 
1964 lectures given at Harvard University. 
Second edition. Lecture Notes in Mathematics,
1500. Springer-Verlag, Berlin, 1992.


\bibitem{Tia} G. Tian,
{\em Smoothness of the universal deformation space of compact
Calabi-Yau manifolds and its Petersson-Weil metric},
in {\bf  Mathematical aspects of string theory}
(San Diego, Calif., 1986), 629--646, Adv. Ser. Math. Phys., 1, 
World Sci. Publishing, Singapore, 1987.

\bibitem{Tod} A.N. Todorov, 
{\em The Weil-Petersson geometry of 
the moduli space of ${\rm SU}(n\geq 3)$
(Calabi-Yau) manifolds. I.},
{\bf Comm. Math. Phys. 126} (1989), no. 2, 325--346. 

\bibitem{Wit} E. Witten, 
{\em Topological sigma models}
{\bf Comm. Math. Phys. 118} (1988), no. 3, 411--449. 

\bibitem{Yau} S.T. Yau ed., 
{\bf Essays on mirror manifolds}, 
 International Press, Hong Kong, 1992.


\end{thebibliography}
\end{document}